\documentclass[12pt, reqno]{amsart}
\usepackage{amsmath, amsthm, amscd, amsfonts, amssymb, graphicx, color}
\usepackage[bookmarksnumbered, colorlinks, plainpages]{hyperref}

\newtheorem{thm}{Theorem}[section]

 \theoremstyle{definition}
 \newtheorem{defn}[thm]{Definition}
 \theoremstyle{remark}

 \numberwithin{equation}{section}
\begin{document}
\title[Pointwise estimate for the Bergman kernel]{Pointwise estimate for the Bergman kernel of the weighted Bergman spaces with exponential type weights}
\author{Sa\"{\i}d Asserda and Amal Hichame}
\address{Ibn tofail university , faculty of sciences, department of mathematics, PO 242 Kenitra Morroco}
\email{asserda-said@univ-ibntofail.ac.ma}

\address{Regional Centre of trades of education and training, Kenitra Morocco}
\email{ amalhichame@yahoo.fr} 
\date{\today}
\begin{abstract}
Let $AL^{2}_{\phi}(\mathbb{D})$ denote the closed subspace of $L^{2}(\mathbb{D},e^{-2\phi}d\lambda)$ consisting of holomorphic functions in the unit disc ${\mathbb D}$. For certain class of   subharmonic funcions $\phi : {\mathbb D}\rightarrow{\mathbb D}$, we prove upper pointwise estimate for the Bergman kernel for $AL^{2}_{\phi}(\mathbb{D})$.
\end{abstract}
\maketitle
\section{Introduction and statement of main result}
\label{}
Let $\mathbb{D}$ be the unit disc in $\mathbb{C}$ and $d\lambda$ be its Lebesgue measure. For a measurable function $\phi : \mathbb{D}\rightarrow\mathbb{D}$, let $L^{2}_{\phi}(\mathbb{D})$ be the Hilbert space of measurable function $f$ on $\mathbb{D}$ such that
$$
\Vert f\Vert_{L^{2}_{\phi}}:=\Bigl(\int_{\mathbb{D}}\vert f\vert^{2}e^{-2\phi}d\lambda\Bigr)^{1\over 2} < \infty
$$
Let $AL^{2}_{\phi}(\mathbb{D})$ be the closed subspace of $L^{2}_{\phi}(\mathbb{D})$ consisting of analytic functions. Let $P$ be the orthogonal projection of $L^{2}_{\phi}(\mathbb{D})$ onto $AL^{2}_{\phi}(\mathbb{D})$ :
$$
Pf(z):=\int_{\mathbb{D}}K(z,w)f(w)e^{-2\phi(w)}d\lambda
$$
where $K$ is the reproducing kernel of $P$.\\ The purpose of this note is to give an upper pointwise estimate of $K$ for some class of subharmonic functions $\phi$ on $\mathbb{D}$ introduced by Oleinik [10] and Oleinik-Perel'man [11].
\begin{defn}
For $\phi \in C^{2}(\mathbb{D})$ and $\Delta\phi > 0$ put $\tau=(\Delta\phi)^{-1/2}$ where $\Delta$ is the Laplace operator. We call $\phi\in \mathcal{OP}(\mathbb{D}) $ if the following conditions holds.\\
(1)\ $\exists\ C_{1} > 0$ such that $\vert\tau(z)-\tau(w)\vert\le C_{1}\vert z-w\vert$ ,\\
(2)\ $\exists\ C_{2} > 0$ such that $\tau(z)\leq C_{2}(1-\vert z\vert)$,\\
(3)\ $\exists\ 0< C_{3} <1$ and $ a > 0$ such that $\tau(w)\leq \tau(z) + C_{3}\vert z-w\vert$ for $ w\notin D(z,a\tau(z))$ where $D(z,a\tau(z))=\{ w\in \mathbb{D}\ ,\ \vert w-z\vert\leq  a\tau(z)\}$.
\end{defn}
Some examples of functions in $\mathcal{OP}(\mathbb{D})$ are as follows :\\
(i) $\phi_{1}(z)=-{A\over 2}\log(1-\vert z\vert^{2}),\ A>0$.\\
(ii) $\phi_{2}(z)={1\over 2}\bigl(-A\log(1-\vert z\vert^{2})+B(1-\vert z\vert^{2})^{-\alpha}\bigr),\ A\geq 0,\ B>0, \alpha > 0$.\\
(iii) $\phi_{1}+h$ and $\phi_{2}+h$ where $\phi_{1}$ and $\phi_{2}$ are as in (i) and (ii) respectively and $h\in C^{2}(\mathbb{D})$ can be any harmonic function on $\mathbb{D}$.\\
For $z,w\in\mathbb{D}$, the distance $d_{\phi}$ induced by the metric $\tau(z)^{-2}dz\otimes d{\bar z}$ is given by
$$
d_{\phi}(z,w)= \inf_{\gamma}\int_{0}^{1}{\vert\gamma^{'}(t)\vert\over\tau(\gamma(t))}dt
$$
where $\gamma$ runs over the piecewise $C^{1}$ curves $\gamma : [0,1]\rightarrow\mathbb{D}$ with $\gamma(0)=z$ and $\gamma(1)=w$. Thanks to condition $(2)$ the metric space $(\mathbb{D},d_{\phi})$ is complete and $d_{\phi}\succeq d_{h}$ where $d_{h}$ is the hyperbolic distance.\\
 Our main result is the following theorem on the off-diagonal decay of the Bergman kernel.
\begin{thm}
Let $\phi\in \mathcal{OP}(\mathbb{D})$ and $K$ be the Bergman kernel for $AL^{2}_{\phi}(\mathbb{D})$. There exist positive constants $C$ and $\sigma$ such that for any $z,w \in \mathbb{D}$
$$
\vert K(z,w)\vert e^{-(\phi(z)+\phi(w))}\leq C{1\over\tau(z)\tau(w)}\exp(-\sigma d_{\phi}(z,w))
$$
\end{thm}
In [4] and [9] M.Christ and J.Marzo-J.Ortega-Cerd\`a obtained a pointwise estimates for the Bergman kernel of the weighted Fock space $\mathcal{F}_{\phi}^{2}(\mathbb{C})$  under the hypothesis that $\Delta\phi$ is a doubling measure. This result was extended to several variables by H.Delin and H.Lindholm in [5] and  [7] under similar hypothesis.\\ In [12], A.P.Schuster and D.Varolin obtained a pointwise estimate for the Bergman Kernel of the weighted Bergman space $AL^{2}(\mathbb{D},e^{-2\phi}(1-\vert z\vert^{2})^{-2}d\lambda)$ under the hypothesis that $\Delta\phi$ is comparable to hyperbolic metric of $\mathbb{D}$ : $$\vert K(z,w)\vert e^{-(\phi(z)+\phi(w))}\leq C\exp(-\sigma d_{h}(z,w))$$
For $\phi\in \mathcal{OP}(\mathbb{D})$ and under the strong condition :  $\forall\ m\geq 1,\ \exists\ b_{m}>0$ and $0<t_{m}<{1\over m}$ such that
$$\tau(w)\leq\tau(z)+t_{m}\vert z-w\vert\ \ \hbox{if}\ \ \vert z-w\vert > b_{m}\tau(z),
$$ H.Arroussi and J.Pau [1] give the following pointwise estimate : for each $k\ge 1$ there exists $C_{k} > 0$ such that
$$\vert K(z,w)\vert e^{-(\phi(z)+\phi(w))}\leq{C_{k}[d_{\tau}(z,w)]^{-k}\over\tau(z)\tau(w)}$$
where $d_{\tau}(w)={\vert z-w\vert\over\min[\tau(z),\tau(w)]}$. A better estimate will be
$$\vert K(z,w)\vert e^{-(\phi(z)+\phi(w))}\leq{C\over\tau(z)\tau(w)}e^{-\sigma d_{\tau}(z,w)}.$$

\section{Proof of theorem 1.2}
\noindent Near the diagonal, by [8,lemma 3.6 ] there exists $\alpha > 0$ sufficiently small such that
$$
\vert K(z,w)\vert\sim\sqrt{K(z,z)}\sqrt{K(w,w)}\sim {e^{\phi(z)+\phi(w)}\over\tau(z)\tau(w)}\quad\hbox{if}\quad \vert z-w\vert\leq\alpha\min[\tau(z),\tau(w)]
$$
Off the diagonal, let $\vert z-w\vert > \alpha\min[\tau(z),\tau(w)]$ and $\beta > 0$ such that $D(z,\beta\tau(z))\cap D(w,\beta\tau(w))=\emptyset$. We may  suppose that $\tau(z)\leq\tau(w)$. Fix a smooth function  $\chi\in C^{\infty}_{0}(\mathbb{D})$ such that \\
- $\hbox{supp}\chi\subset D(w,\beta\tau(w))$,\\
- $0\leq\chi\le 1$, $\chi=1$ in $D(w,{\beta\over 2}\tau(w))$ and \\ 
- $\vert\bar\partial\chi\vert^{2}\preceq\chi\tau(w)^{-2}$.\\ Since $\phi\in\mathcal{OP}(\mathbb{D})$, by [10,lemma 1 and 2] the following mean inequality holds
\begin{eqnarray*}
(*)\qquad\vert K(w,z))\vert^{2}e^{-2\phi(w)}&\preceq&{1\over\tau(w)^{2}}\int_{D(w,{\beta\over 2}\tau(w))}\chi(\zeta)\vert K(\zeta,z))\vert^{2}e^{-2\phi(\zeta)}d\lambda(\zeta)\\
&\preceq&{1\over\tau(w)^{2}}\Vert K(.,z)\Vert^{2}_{L^{2}(\chi e^{-2\phi}d\lambda)}
\end{eqnarray*}
Hence $$\Vert K(.,z)\Vert_{L^{2}(\chi e^{-\phi})}=\sup_{f}\vert<f,K(.,z)>_{L^{2}(\chi e^{-2\phi}d\lambda)}\vert$$ where $f$ is holomorphic in $D(w,\beta\tau(w))$ with $\Vert f\Vert_{L^{2}(\chi e^{-2\phi}d\lambda)}=1$. Since $P_{\phi}(f\chi)(z)=<f,K(.,z)>_{L^{2}(\chi e^{-2\phi}d\lambda)}$ and  that $u_{f}=f\chi-P_{\phi}(f\chi)$ is the minimal solution in $L^{2}(\mathbb{D}, e^{-2\phi}d\lambda)$ of $\bar\partial u=f\bar\partial\chi$, and from the fact  $\chi(z)=0$, we have
$$
\vert<f,K(.,z)>_{L^{2}(\chi e^{-2\phi}d\lambda)}\vert=\vert P_{\phi}(f\chi)(z)\vert=\vert u_{f}(z)\vert
$$
Since $D(z,\beta\tau(z))\cap D(w,\beta\tau(w))=\emptyset$, the function $u_{f}$ is holomorphic in $D(z,\nu\tau(z))$ for some $\nu > 0$. By the mean value inequality
\begin{eqnarray*}
\vert u_{f}(z)\vert^{2}e^{-2\phi(z)}&\preceq&{1\over\tau(z)^{2}}\int_{D(z,\nu\tau(z))}\vert u_{f}(\zeta)\vert^{2}e^{-2\phi(\zeta)}d\lambda\\
&\preceq&{1\over\tau(z)^{2}}\int_{D(z,\nu\tau(z))}e^{-\epsilon{\vert\zeta-z\vert\over\nu\tau(z)}}\vert u_{f}(\zeta)\vert^{2}e^{-2\phi(\zeta)}d\lambda
\end{eqnarray*}
 Since the linear curve $\gamma(t)=(1-t)z+t\zeta$ lies in $D(z,\nu\tau(z))$ and $\tau(\gamma(t))\sim\tau(z)$, we have $d_{\phi}(\zeta,z)\le C{\vert \zeta-z\vert\over\tau(z)}$ for $\zeta\in D(z,\nu\tau(z))$. Hence
\begin{eqnarray*}
\vert u_{f}(z)\vert^{2}e^{-2\phi(z)}&\preceq &{1\over\tau(z)^{2}}\int_{D(z,\nu\tau(z))}e^{-C\epsilon d_{\phi}(\zeta,z)}\vert u_{f}(\zeta)\vert^{2}e^{-2\phi(\zeta)}d\lambda\\
&\preceq &{1\over\tau(z)^{2}}\int_{\mathbb{D}}e^{-C\epsilon d_{\phi}(\zeta,z)}\vert u_{f}(\zeta)\vert^{2}e^{-2\phi(\zeta)}d\lambda
\end{eqnarray*}
The function $\zeta\rightarrow d_{\phi}(\zeta,z)$ is smooth on $\mathbb{D}\setminus\hbox{Cut}(z)\cup\{z\}$ where $\hbox{Cut}(z)$ is the cut locus : the set of all cut points of $z$ along all geodesics that start from $z$. To get a smooth Lipschitz approximation of $d_{\phi}$, we recall the following result of Greene-Wu [6] ( see also [2]).
\begin{thm}
Let $M$ be a complete Riemannian manifold, let $ h : M\rightarrow\mathbb{R}$ be a Lipschitz function, let $\eta : M\rightarrow ]0,+\infty[$ be a continuous function, and $r$ a positive number. Then there exist a smooth Lipschitz function $ g : M\rightarrow\mathbb{R}$ such that $\vert h(x)-g(x)\vert\leq\eta(x)$ for every $x\in M$, and $\hbox{Lip}(g)\leq\hbox{Lip}(h)+r$.
\end{thm}
We use this result with $h(\zeta)=d_{\phi}(\zeta,z),\ \eta=1$ and $r=1$. We have $d_{\phi}(\zeta,z)\prec g_{z}(\zeta)\prec d_{\phi}(\zeta,z)$ and $\tau(\zeta)\vert dg_{z}(\zeta)\vert\leq 2$. Hence
$$
\vert u_{f}(z)\vert^{2}e^{-2\phi(z)}\preceq{1\over\tau(z)^{2}}\int_{\mathbb{D}}e^{-C\epsilon g_{z}(\zeta)}\vert u_{f}(\zeta)\vert^{2}e^{-2\phi(\zeta)}d\lambda
$$
By Berndtsson-Delin's improved $L^{2}$ estimates of for the minimal solution of $\bar\partial\ $ in $L^{2}(\mathbb{D},e^{-2\phi}d\lambda)$ [3][5] , we have :
$$
\int_{\mathbb{D}}e^{-C\epsilon g_{z}(\zeta)}\vert u_{f}(\zeta)\vert^{2}e^{-2\phi(\zeta)}d\lambda\preceq\int_{\mathbb{D}}e^{-C\epsilon g_{z}(\zeta)}\vert \bar\partial \chi(\zeta)\vert^{2}\vert f(\zeta)\vert^{2}\tau(\zeta)^{2}e^{-2\phi(\zeta)}d\lambda
$$
provided that $\tau\vert\partial\omega_{\epsilon}\vert\leq\mu\omega_{\epsilon}$ with $\mu < \sqrt{2}$ where $\omega_{\epsilon}(\zeta)=e^{-C\epsilon g_{z}(\zeta)}$. If we choose $\epsilon$ small enough so that $\mu=2C\epsilon<\sqrt{2}$ then $\tau\vert\partial\omega_{\epsilon}\vert=C\epsilon\tau\vert\partial g_{z}\vert\omega_{\epsilon}\leq \mu\omega_{\epsilon}$. Thus
$$
\vert u_{f}(z)\vert^{2}e^{-2\phi(z)}\preceq{1\over\tau(z)^{2}}\int_{D(w,\beta\tau(w))}e^{-C\epsilon d_{\phi}(\zeta,z)}\chi(\zeta)\vert f\vert^{2}e^{-2\phi(\zeta)}d\lambda
$$
where for the last term we use $\tau(\zeta)\sim\tau(w)$. Since $\zeta\in D(w,\beta\tau(w))$ we have
$$
d_{\phi}(\zeta,z)\geq d_{\phi}(z,w)-d_{\phi}(w,\zeta)\succeq d_{\phi}(z,w)-{\vert\zeta-w\vert\over\beta\tau(w)}\succeq d_{\phi}(z,w)
$$
and thanks to $(*)$, we conclude
$$
\vert K(z,w)\vert e^{-(\phi(w)+\phi(z))}\leq{C\over\tau(z)\tau(w)}e^{-\sigma d_{\phi}(z,w)}.
$$

\vskip 20 pt

\end{document}